\newcount\notenumber

\def\note{\advance\notenumber by 1
\footnote{$^{(\the\notenumber)}$}}

\def\A{{\bf A}}
\def\Z{{\bf Z}}
\def\Q{{\bf Q}}
\def\C{{\bf C}}
\def\O{{\cal O}}
\def\N{{\bf N}}
\def\Pr{{\bf P}}

\def\ord{{\rm ord}}

\font\title=cmr10 scaled 1200

%30-4-2002

\centerline{\title On integral points on surfaces}\bigskip

\centerline{P. Corvaja\hskip 2truecm U. Zannier}\bigskip

\noindent{\bf Abstract.}\ We  study the integral points
on surfaces by means of a new method, relying on the Schmidt Subspace
Theorem. This
method was recently introduced in [CZ] for the case of curves, leading to a new
proof of Siegel's celebrated theorem. Here, under certain conditions
involving the
intersection matrix of  the divisors at infinity, we shall conclude that
the integral
points on a surface all lie on a curve. We shall also give several examples and
applications. One of them concerns  curves, with a study of the
integral points  defined over a variable quadratic field; for instance we
shall show that an affine curve with at least five points at infinity has
at most finitely many such integral points.\bigskip

\noindent{\bf \S 0. Introduction and statements.} In the recent paper [CZ]
a new
method was introduced in connection with the integral points on an
algebraic curve;
this led to a novel proof of Siegel's celebrated theorem, based on the Schmidt
Subspace Theorem  and entirely avoiding any recourse to abelian varieties
and their
arithmetic.
 Apart from this methodological point, we observed (see the Remark in [CZ])
that the
approach was sometimes capable of quantitative improvements on the
classical one, and
we also alluded to the possibility of extensions to higher dimensional
varieties. The
present paper represents precisely a first step in that direction, with an
analysis
of the case of surfaces. \medskip

The crux of the arguments in [CZ] appeared through the special case of Siegel's
Theorem  when the affine curve
misses at least {\it three} points with respect to its projective closure; this
condition alone implies the finiteness of the set of integral points. That
case was
studied by embedding the curve  in a  space of
 large dimension and  by constructing hyperplanes with high order contact
with the curve  at some point at infinity; finally one exploited the
diophantine
approximation through the Schmidt Theorem rather than through the Roth's
one, as in
the usual approach. Correspondingly, here we shall work with (nonsingular)
affine
surfaces missing at least {\it four} divisors; but now, unlike the case of
curves, we
shall need additional assumptions on the divisors, expressed through their
intersection matrix. These  conditions appear  naturally when using the
  Riemann-Roch Theorem to embed the
surface in a suitable space  and to construct functions with zeros of large
order
along a prescribed divisor in the set,    allowing an  application of the
Subspace
Theorem.

The result of this approach is  the Main Theorem below. Its
assumptions appear somewhat technical, so we have preferred to start with its
corollary Theorem 1 below; this is sufficient for some applications, such as to
Corollary 1, which concerns the {\it quadratic integral} points on a curve.
As a
kind of ``test" for the Main Theorem, we shall see how it
immediately implies Siegel's theorem on curves (Ex. 1.5). Still other
applications of
the method may be obtained looking at varieties defined in $\A^m$ by one
equation
$f_1\cdots f_r=g$, where $f_i,g$ are polynomials and $\deg g$ is ``small".
(A special
case arises with ``norm form equations", treated by Schmidt in full
generality; see
[S1].) However  in general  the variety has singularities at infinity, so, even in the case
of surfaces, 
the Main
Theorem cannot be applied directly to such equation; this is why we postpone such analysis to a
separate paper.  \medskip

In the sequel we let $\tilde X$ denote a geometrically irreducible non-singular
projective surface defined over a number field $k$. We also let $S$ be  a
finite set
of places of $k$, including the archimedean ones, denoting as usual
$\O_S=\{\alpha\in k: |\alpha|_v\le 1\quad \hbox{for all $v\not\in S$}\}$.

We view the {\it
$S$-integral points} in the classical way; namely, letting  $X$ be an affine
Zariski-open subset of $\tilde X$ (defined over $k$), embedded in $\A^m$,
say, we
define an $S$-integral point $P\in X(\O_S)$ as a point whose coordinates lie in
$\O_S$. For our purposes, this is equivalent with the more modern
definitions given e.g. in [Se1] or [V].\medskip

\noindent{\bf Theorem 1.}\enspace {\it Let $\tilde X$ be a surface as
above, and let $X\subset\tilde X$ be an affine open subset. 
Assume that $\tilde X\setminus X=D_1\cup\ldots \cup D_r$, 
where the $D_i$ are distinct irreducible divisors such
that no three of them share a common point. 
Assume also that $r\ge 4$ and that there exist positive integers
$p_1,\ldots ,p_r,c$,  such that $p_ip_j(D_i.D_j)=c$ for all pairs $i,j$.

Then there exists a curve on $X$ containing all the  $S$-integral points in
$X(k)$.}\medskip

Below  we shall note that one cannot remove the
condition on the $(D_i.D_j)$ (see Ex. 1.1).\medskip

An application of Theorem 1 (which does not follow from the mentioned
results by Vojta), concerns the  points on a curve which are integral and
defined over a field of degree at most $2$  over $k$; 
we insist that here we do not view this field as being fixed, but varying
with the point.  This situation (actually for fields of any given degree in
place of $2$)  has been studied in the context of
rational points, via the former Mordell-Lang conjecture, now proved by
Faltings; see e.g. [HSi, pp. 439-443] for  an account of some results and
several references. For instance, in the quadratic case it follows from
rather general results by D. Abramovitch and J. Harris (see [HSi, Thms F121,
F125(i)]) that {\it if a curve has infinitely many  points
rational over a quadratic extension of $k$, then it admits a map of degree
$\le 2$ either to $\Pr^1$ or to an elliptic curve}.

For integral points we may obtain   without appealing to
Mordell-Lang  a result in the same vein, which however seems not to derive
directly from the rational case, at least when the genus is $\le 2$. 
(In fact, in that case, Mordell-Lang  as applied in [HSi] gives no
information at all.)   This result will be
proved by applying Theorem 1 to the symmetric product of a curve with
itself. We state it as a corollary, where we use the terminology 
{\it quadratic (over $k$) $S$-integral  point} 
to mean  a point defined over a quadratic extension
of $k$, which is integral at all places of $\overline\Q$
except possibly  those lying above  $S$. \medskip

\noindent{\bf Corollary 1.}\enspace {\it Let $\tilde C$ be a geometrically
irreducible projective curve and let $C=\tilde C\setminus \{A_1,\ldots ,A_r\}$ 
be an affine subset, where the $A_i$ are distinct
points in $\tilde C(k)$. Then

(i) If $r\ge 5$, $C$ contains only finitely many quadratic (over $k$)
$S$-integral
points.

(ii) If $r\ge 4$, there exists a finite set of rational maps $\psi :\tilde
C\rightarrow \Pr^1$ of degree $2$ such that all but finitely many of the
quadratic $S$-integral points on $C$  are sent to  $\Pr^1(k)$ by some of the
mentioned maps.}\medskip

In the next section we shall see that the result is in  a sense
best-possible (see Ex. 1.2-1.3), 
and we shall briefly discuss possible extensions.  We shall also
state an ``Addendum" which provides further information on the maps  in
(ii). \medskip

%%%%%%%%%%%%%%%%%%%%%%%%%%%%%%%%%%%%%%%%%%%%%%%%%%%%%%%%%%%%%%%%%%%%%%%%%%%%%%

As mentioned  earlier,  we have postponed the statement of our
main result (which implies Theorem 1), because of its   somewhat
involved formulation. Here it is: \medskip

\noindent{\bf Main Theorem.}\enspace {\it Let $\tilde X$ be a surface as
above, and
let $X\subset\tilde X$ be an affine open subset. Assume that $\tilde X\setminus
X=D_1\cup\ldots \cup D_r$, $r\ge 2$, where the $D_i$ are distinct
irreducible divisors
with the following properties:\smallskip

\item{(i)} -  No three of the $D_i$ share a common point.

\item{(ii)} - There exist positive integers $p_1,\ldots ,p_r$ such that, 
putting $D:=p_1D_1+\ldots +p_rD_r$, $D$ is ample and the following holds. 
Defining $\xi_i$, for $i=1,\ldots,r$, as the minimal positive 
 solution of the equation  $D_i^2\xi^2-2(D.D_i)\xi+D^2=0$ 
($\xi_i$ exists; see \S 2),  we have   the inequality
$$
2D^2\xi_i>(D.D_i)\xi_i^2+3D^2p_i.
$$ \smallskip

Then there exists a curve on $X$ containing all the  $S$-integral points in
$X(k)$.}\bigskip

Our proofs, though not effective in the sense of
leading  to explicit equations for the relevant curve, allow in principle
 quantitative conclusions such as an explicit estimation of the degree of
the curve.
Also,  the bounds may be obtained to be rather uniform with respect to the
field $k$; one may use  results due to Schlickewei,
Evertse  (as for instance in the Remark in [CZ], p. 271) or more recent
estimates by Evertse and Ferretti [EF]; 
this last paper uses the proof-approach to the Subspace Theorem due
to Faltings and W\"ustholz [FW], through the product theorem [F]. However
here we shall not pursue in this direction.

\bigskip

%%%%%%%%%%%%%%%%%%%%%%%%%%%%%%%%%%%%%%%%%%%%%%%%%%%%%%%%%
%%%%%%%%%%%%%%  Dimostrazione Teor. generale
%%%%%%%%%%%%%%%%%%%%%%%%%%%%%%%%%%%%%%%%%%%%%%%%%%%%%%%%%

\noindent{\bf \S 1. Remarks and examples.} In this section we collect several
observations on the previous statements. Concerning Theorem 1, we start by
pointing
out that the condition on the $(D_i.D_j)$ cannot be removed.\medskip

\noindent{\bf Example 1.1.} Let $\tilde X=\Pr^1\times \Pr^1$ and let
$D_1,\ldots ,D_4$ 
be the divisors $\{0\}\times\Pr^1$, $\{\infty\}\times\Pr^1$, $\Pr^1\times
\{0\}$ and $\Pr^1\times\{\infty\}$ in some order. Then, defining $X:=\tilde
X\setminus (\cup_{i=1}^4D_i)$, we see that $X$ is isomorphic to the product
of the affine line minus one point with itself. 
Therefore the integral points on $X$ are
(for suitable $k,S$) Zariski dense on $X$. (On the contrary, Theorem 1
easily implies that the integral points on $\Pr^2$ 
minus four divisors in general position are not Zariski dense, 
a well-known fact.)\medskip

Theorem 1 intersects results  due to Vojta (also obtained through the Subspace
Theorem); see e.g. [V, Thms. 2.4.1, 2.4.6]; roughly speaking, these statements
predict that the integral points are not Zariski dense, provided there are
sufficiently many components at infinity. However, they  do not lead
directly to the general case of Theorem 1, as happens  for instance when
Pic$^o(\tilde X)$  is not trivial (like in the proof of Corollary 1), or
when  the rank of the subgroup generated by the $D_i$ in  
 $NS(\tilde X)$ is greater than $1$.

The conditions on the number   of divisors $D_i$ and on the $(D_i.D_j)$
which appear in the Main Theorem (and in Theorem 1)  come
naturally from our method. One may ask how these assumptions  fit with
celebrated conjectures on integral points (see [HSi, Ch. F]). 
We do not have any definite view here; 
we just recall Lang's point of view, expressed in [L, p. 225-226];
namely, on the one hand Lang's Conjecture 5.1, [L, p. 225], 
 predicts at most finitely many integral points on hyperbolic varieties; 
on the other hand, it is ``a general idea"
that taking out a sufficiently large number of divisors 
(or a divisor of large degree)
from a projective variety produces a hyperbolic space. 
Lang interprets in this way also the  results by Vojta alluded to above. 
\note{Both our method and Vojta's do not work at all by removing 
a single divisor (but see Ex. 1.4 below).} \medskip

We now turn to Corollary 1, noting  that in
some sense its conclusions   are best-possible.\medskip

\noindent{\bf Example 1.2.}  Let  a rational map $\psi:\tilde
C\rightarrow\Pr^1$ of
degree 2 be given. We construct an
affine subset $C\subset\tilde C$ with four missing points and infinitely many
quadratic integral points. Let $B_1,B_2$ be distinct points in $\Pr^1(k)$
and define
$Y:=\Pr^1\setminus \{B_1,B_2\}$. Lifting $B_1,B_2$ by $\psi$ gives in
general four
points $A_1,\ldots ,A_4\in \tilde C$. Define then  $C=\tilde
C\setminus\{A_1,\ldots
,A_4\}$. Then $\psi$ can be seen as a finite morphism from $C$ to $Y$.
Lifting (the
possibly) infinitely many integral points in $Y(\O_S)$ by $\psi$ produces then
infinitely many quadratic   $S'$-integer  points on $C$ (for a suitable finite
$S'\supset S$).\medskip

 Concrete examples are obtained e.g. with  the classical space curves
given by two simultaneous  Pell equations, such as e.g.  $t^2-2v^2=1$,
$u^2-3v^2=1$.
We now have an affine subset of an elliptic curve, with four points at
infinity. We
can obtain infinitely many quadratic integral points by solving in $\Z$
e.g.  the
first Pell  equation, and then defining $u=\sqrt{3v^2+1}$; or we may solve
the second
equation and then put $t=\sqrt{2v^2+1}$; or we may also  solve
$3t^2-2u^2=1$ and then
let $v=\sqrt{t^2+1\over 2}$. (This is the  construction of Example 1.2  for
the three natural projections.)

It is actually possible to show through Corollary 1 that all but finitely many
quadratic integral points arise in this way.
\note{On the contrary, the quadratic rational points  cannot be likewise
described;
we can obtain them as inverse images from $\Pr^1(k)$ under any map of degree  2
defined over $k$, and it is easy to see that  in general no finite set  of
such maps
is sufficient to obtain almost all the points in question.}  We in fact have an
additional property for the  relevant maps in conclusion (ii), namely:\medskip

\noindent{\bf Addendum to Corollary 1.}\enspace {\it Assume that $\psi$ is
a quadratic map as in
(ii) and that it sends to $\Pr^1(k)$ an infinity of the integral points in
question.
Then the set $\psi (\{A_1,\ldots ,A_4\})$ has two points. In particular, we
have a linear-equivalence relation $\sum_{i=1}^4\epsilon_i(A_i)\sim 0$ on
Div($\tilde
C$), where the $\epsilon_i\in \{\pm 1\}$ have zero sum.}\medskip

 When such a $\psi$ exists, the two
relevant values of it can  be sent to two prescribed points in $\Pr^1(k)$ by
means of an automorphism of $\Pr^1$; in practice, the choice of the maps
$\psi$ then
reduces to splitting the four points at infinity in two pairs having equal
sum in the
Jacobian of $\tilde C$; this can be done in at most three
ways, as in the example with the Pell equations. The simple proof for the {\it
Addendum}  will be given after the one for the corollary. This conclusion
of course
allows one to compute the relevant maps and to  parametrize  all but finitely
many quadratic integral points
on an affine curve with four points at infinity.  \medskip

 Concerning again Cor. 1 (ii),  we   now observe that
   ``$r\ge 4$" cannot be substituted with $r\ge 3$. \medskip

\noindent{\bf Example 1.3.}\ Let $C=\Pr^1\setminus\{-1,0,\infty\}$,
realized  with the
plane equation $X(X+1)Y=1$. Let $r,s$
run through
the $S$-units in $k$ and define $a={r-s+1\over 2}$,  $\Delta =a^2-r$. Then
the points
given by $x=a+\sqrt \Delta$, $y={x'(x'+1)\over rs}$, where
$x'=a-\sqrt\Delta$, are quadratic $S$-integral on $C$. It is possible to show  that
they  cannot all be mapped to $k$ by one at least of a finite number of 
quadratic  maps. \medskip

It is also possible to show that for the affine elliptic curve $E: Y^2=X^3-2$, the
quadratic integral points (over $\Z$) cannot be all described like in (ii) of
Corollary 1.

Note that $E$ has only one point at infinity.  Probably similar examples
cannot  be constructed with more points at infinity; namely, (ii) is
unlikely to
be  best-possible also for curves of   genus $g\ge 1$, in the sense that the
condition $r\ge 4$ may be then probably relaxed. In fact, a conjecture of
Lang and
Vojta (see [HSi, Conj. F.5.3.6, p. 486]) predicts that {\it if  $X=\tilde
X\setminus
D$ is an affine variety with $K_X+D$ almost ample} (i.e. ``big") {\it and
$D$ with
normal crossings, the integral points all lie  on a proper subvariety}.
Now, in the
proof of our corollary we work with $\tilde X$ equal to $\tilde C^{(2)}$, the
two-fold  symmetric power of $\tilde C$, and with $D$ equal to the image in
$\tilde C^{(2)}$  of $\sum_{i=1}^rA_i\times \tilde C$. 
It is then easily checked  that $K_X+D$ is (almost) ample precisely when $g=0$
and $r\ge 4$,  or $g=1$ and $r\ge 2$ or $g\ge 2$ and $r\ge 1$. 
In other words, the Lang-Vojta conjecture essentially
predicts that counterexamples sharper than those given here
may not be found.

 To prove this, one might try to proceed like in the deduction of   
Siegel's Theorem from the special case of three points at infinity. Namely, one may
then use  unramified covers, as in [CZ], with the purpose of  increasing  the number of
points at infinity. (One also uses [V, Thm. 1.4.11], essentially the Chevalley-Weil
Theorem, to show that lifting the integral  points does not produce infinite degree
extensions.) 

In the case of the present Corollary 1 a similar strategy does not help.
In fact,  the structure of the fundamental group of $\tilde C^{(2)}$
\note{Angelo Vistoli has pointed out to us that it is the abelianized of 
$\pi_1(\tilde C)$.} prevents the number of components  
of a divisor to increase by pull-back on a cover.
 However
there exist nontrivial instances beyond the case of curves, and showing one of them 
 is our purpose in including this further result, namely:\medskip

\noindent{\bf Example 1.4.}\enspace {\it Let $A$ be an abelian variety of
dimension $2$, let $\pi:A\rightarrow A$ be an isogeny of degree $\ge 4$
 and let $E$ be an ample irreducible divisor on $A$. We suppose that for
$\sigma\in\ker\pi$    no three of the divisors
$E+\sigma$ intersect. Then there are at most finitely many  $S$-integral points in 
$(A\setminus \pi(E))(k)$.}\medskip

We remark that this is an extremely special case of a former conjecture by
Lang, proved by Faltings [F, Cor. to Thm. 2]: {\it every affine subset of
an abelian
variety has at most finitely many integral points}. 
\medskip

We just sketch a proof.  Note now that $\pi(E)$ is an irreducible divisor, so Theorem
1   cannot be applied directly. Consider $D:=\pi^*(\pi(E))$;  since  $\pi$ has degree
$\ge 4$, we see  that $D$ is the sum of $r:=\deg\pi \ge 4$ irreducible divisors
satisfying the assumptions for Theorem 1, with $p_i=1$ for $i=1,\ldots ,r$.

Let now $\Sigma$ be an infinite set of $S$-integral points in $Y(k)$, where
$Y=A\setminus \pi(E)$. By [V, Thm. 1.4.11],  $\pi^{-1}(\Sigma)$ is
a set of  $S'$-integral points on $X(k')$, where $X=A\setminus D$, for some
number field $k'$ and some finite set $S'$ of places of $k'$. 
By Theorem 1 applied to $X$ we
easily deduce the conclusion, since there are no curves of genus zero on an abelian
variety ([HSi, Ex. A74(b)]).///\medskip

We conclude this section by showing how the Main Theorem leads directly to
Siegel's
Theorem for the case of at least three points at infinity.  (As  remarked
above, one recovers the full result by taking, when genus$(C)>0$, an
unramified cover of degree $\ge 3$ and applying  the special case and
 [V, Thm. 1.4.11].)\medskip

\noindent{\bf Example 1.5.} We prove: {\it Let $\tilde C$ be a projective
curve and
$C=\tilde C\setminus\{A_1,\ldots ,A_s\}$, $s\ge 3$ an affine subset. Then
there are at most finitely many $S$-integral points on $C$.} This special case of
Siegel's Theorem appears as Theorem 1 in [CZ]. We now show how this follows
at once from the Main Theorem. 
 First, it is standard that one can reduce to nonsingular
curves. We then let $\tilde X=\tilde C\times\tilde C$ and $X=C\times C$. 
Then $\tilde X\setminus X$ is the union of $2s$ divisors $D_i$ of the form
$A_i\times\tilde C$ or $\tilde C\times A_i$, 
which will be referred to as of the  {\it first} or {\it second} type
respectively. Plainly, the intersection product $(D_i.D_j)$ will be $0$ or $1$
according as $D_i,D_j$ are of equal or different types. We put in the Main
Theorem $r=2s$, $p_1=\ldots =p_r=1$. All the hypotheses are verified except
possibly (ii).
To verify  (ii), note that $(D_i.D_i)=0$,  $(D.D_i)=s$, $D^2=2s^2$. Therefore
$\xi_i=s$ and we have to prove that 
$4s^3>s^3+6s^2$ which is true precisely when  $s>2$.

We conclude that the $S$-integral points on $C\times C$ are not Zariski dense, whence
the assertion.

\bigskip

\noindent {\bf \S 2. Tools from intersection theory on surfaces}. We shall
now recall
a few simple facts from the theory of surfaces, useful for the proof of
Main Theorem.
These include a version of the Riemann-Roch theorem and involve intersection
products. (See e.g. [H, Ch. V] for the basic theory.)\medskip

Let $\tilde{X}$ be a projective smooth algebraic surface defined
over the complex number field $\C$. We will follow the
notations of [B] (especially Chapter 1), which are rather
standard. For a divisor $D$ on $\tilde{X}$ and an integer
$i=0,1,2$, we denote by
$h^i(D)$ the dimension of the vector space ${\rm H}^i(\tilde{X},\O(D))$.
We shall make essential use of the following  asymptotic version
of the Riemann-Roch  theorem:
\medskip

\noindent {\bf Lemma 2.1}
{\it Let $D$ be an ample divisor on $\tilde{X}$.
Then for positive integers $N$ we have
$$
h^0(ND)={N^2D^2\over 2}+O(N).
$$
}

\noindent {\it Proof}. The classical Riemann-Roch theorem
(see e.g. Th\'eor\`eme I.12 of [B] and the following Remarque
I.13) gives
$$
h^0(ND)={1\over 2} (ND)^2-{1\over 2}(ND.K)+\chi(\O_X)+h^1(ND)-h^0(K-ND),
$$
where $K$ is a canonical divisor of $\tilde{X}$.
The first term is precisely $N^2D^2/2$. Concerning the other terms, observe
that:
$h^1(ND)$ and $h^0(K-ND)$
vanish for large $N$; $\chi(\O_X)$ is  constant; the intersection
product $(ND.K)$ is linear in $N$. The result  then follows. ///
\medskip

We will need an estimate for the dimension of the linear
space of sections of ${\rm H}^0(X,\O(ND))$ which have a zero
of given order on a fixed (effective) curve $C$. We begin with a lemma.
\medskip

\noindent {\bf Lemma 2.2}. {\it Let $D$ be a divisor, $C$ a
  curve on $\tilde{X}$; then
$$
h^0(D)-h^0(D-C)\leq\max\{0,1+(D.C)\}.
$$
}

\noindent {\it Proof}. In proving the inequality we may replace $D$ with any
divisor linearly equivalent to it. In particular, we may assume that $|D|$
does not
contain  any possible singularity of $C$.

  Let us then recall that for
every sheaf  $\cal L$ the exact sequence $$ 0\rightarrow {\cal
L}(-C)\rightarrow {\cal
L}\rightarrow {\cal L}|C\rightarrow 0
$$
gives an exact sequence in cohomology
$$
0\rightarrow
{\rm H}^0(\tilde{X},{\cal L}(-C))\rightarrow
{\rm H}^0(\tilde{X},{\cal L})\rightarrow
{\rm H}^0(C,{\cal L}|C)\rightarrow \ldots
$$
from which we get
$$
\dim({\rm H}^0(\tilde{X},{\cal L})/{\rm H}^0(\tilde{X},{\cal L}(-C)))
\leq
\dim {\rm H}^0(C,{\cal L}|C).
$$
Applying this inequality with ${\cal L}=\O(D)$
we get
$$
h^0(D)-h^0(D-C)\leq \dim {\rm H}^0(C,\O(D)|C).
$$
The sheaf $\O(D)|C$ is an invertible sheaf of
degree $(D.C)$ on the   complete curve $C$.
(See [B, Lemme 1.6], where $C$ is assumed to be smooth; this makes no
difference
because of our opening assumption on $|D|$.) We can then bound the right
term by $\max\{0, 1+(D.C)\}$ as wanted.///
\medskip

\noindent{\bf Lemma 2.3}. {\it Let $D$ be an ample effective divisor on
$\tilde{X}$, $C$ be an irreducible component of $D$.
For  positive integers $N$ and $j$ we have that either
${\rm H}^0(\tilde{X},\O(ND-jC))=\{0\}$ or
$$
0\leq h^0(ND-jC)-h^0(ND-(j+1)C)\leq N(D.C)-jC^2+1.
$$
 }

\noindent {\it Proof}. Suppose first that $(ND-jC.C)\geq 0$. Then
 Lemma 2.2 applied with $ND-jC$ instead of $D$ gives
what we want. If otherwise $ND-jC$ has negative intersection with
the effective curve $C$ then $\O(ND-jC)$ has no regular sections. In fact,
assume the contrary. Then  there would exist an effective divisor $E$
linearly equivalent to $ND-jC$, whence $E.C=(ND-jC.C)<0$. But $E.C$ must be
$\ge 0$. (In fact,  since $E$ is effective we may write $E=E_1+rC$, where $E_1$ is
effective and does not contain $C$ and where $r\ge 0$. Then $E.C=E_1.C+rC^2$,
whence the claim, in view of  $C^2>0$, which in turn follows from
$(ND-jC.C)<0$).
This contradiction concludes the proof./// \bigskip

\noindent{\bf Lemma 2.4} {\it Let $D$ be an ample divisor,
$C$ be an effective curve. Then
$$
D^2C^2\leq (D.C)^2.
$$}

\noindent {\it Proof}. This is in fact well known (see e.g. [H, Ch. V, Ex.
1.9]).
We give however a short proof for completeness. The inequality is non
trivial only in the case $C^2>0$. 
Assume this holds. Then if we had $D^2C^2> (D.C)^2$, the
intersection form on the rank two group generated by $D$ and $C$ in ${\rm
Pic}(\tilde{X})$ would be positive definite, which contradicts the Hodge index
theorem [H, Ch. V, Thm. 1.9].///
\bigskip

%%%%%%%%%%%%%%%%%%%%%%%%%%%%%%%%%%%%%%%%%%%%%%%%
 
%%%%%%%%%%%%%%%%%%%%%%%%%%%%%%%%%%%%%%%%%%%%%%%

\noindent{\bf \S 3. Proofs.} We shall begin with the proof of Main Theorem,
actually anticipating a few words on the strategy. Then we shall deduce
Theorem 1
from the Main Theorem. In turn, Theorem 1 shall be employed for the proof of
Corollary 1. \bigskip

\noindent{\it Proof of Main Theorem.} We begin with a brief sketch of our
strategy, assuming for simplicity that $S$ consists of just   
one (archimedean) absolute value.
In the case treated in [CZ], of an affine curve $C$ with   missing points
$A_1,\ldots ,A_r$, $r\ge 3$, we first embed $C$ in a high dimensional space
by means of a basis for the space $V$ of regular functions on $C$ with at most
poles of order $N$ at the given points. 
Then, going to an infinite subsequence  $\{P_i\}$ of the integral
points on $C$, we may assume that $P_i\rightarrow A$, where $A$ is some $A_i$.
Linear algebra now gives functions in $V$ vanishing at $A$ with orders 
$\ge -N$, $\ge -N+1,\ldots ,\ge -N+d$, where $d=\dim V$. 
Such functions may be viewed as
linear forms in the previous basis and these vanishings imply that the
product of these functions evaluated at the $P_i$ is small. 
Then the Subspace Theorem (recalled below) applies.

The principles are similar in the present case of surfaces, the role of the
points $A_i$ being now played by the divisors $D_i$. However one has to deal
with several new technical difficulties. For instance, the construction of the
functions with large order zeros is no longer automatic and the quantification
involves now intersection indices. Moreover,   additional complications
appear  when
the integral points converge simultaneously to two divisors in the set,
i.e. to some
intersection point (this is ``Case C" of the proof below). 
%This affects the final
%estimates, and it turns out that in this case it is sometimes convenient to
%construct not functions with large vanishing along a whole divisor, but rather
%construct functions with ``large index" at the given point.  \bigskip

Now we go on with the details. We shall assume throughout that each of the
divisors $D_i$ is defined over $k$. Also, we assume that each valuation
$|\cdot|_v$
is  normalized so that if $v|p$, then $|p|_v=p^{-{[k_v:\Q_p]\over
[k:\Q]}}$,  where
$k_v$ is the completion of $k$ at $v$, and similarly for archimedean $v$.
As usual,
for a point $(x_1:\ldots :x_d)\in{\bf P}^{d-1}(k)$, ($d\ge 2$), we define the
projective
 height as  $H(x_1:\ldots :x_d)=\prod_v\max(|x_1|_v,\ldots
,|x_d|_v)$.\medskip

The theorem will follow if we prove that {\it for every
infinite sequence of integral points on $X$, there exists a curve
defined over $k$ containing an infinite subsequence}. In fact,
 arrange all the curves on $X$ defined over $k$ in
a sequence $C_1,C_2,....$ Now, if the conclusion of the theorem  is
not true, we may find for each $n$ an integral point $P_n$
on $X$ outside $C_1\cup C_2\cup\ldots \cup C_n$. But then no given
curve $C_m$ can contain infinitely many of the points $P_i$.\medskip

Let then $\{P_i\}_{i\in \N}$ be an infinite sequence of pairwise
distinct integral points on  $X$. By the observation just made, we may
restrict our attention  to any infinite subsequence, and thus we may
assume in particular that for each valuation $v\in S$  the $P_i$
converge $v$-adically to a point $P^v\in \tilde X(k_v)$. \medskip

We recall that $D_i$, $i=1,\ldots ,r$, are certain irreducible divisors
on $\tilde X$, and that we put $D=\sum_{i=1}^rp_iD_i$, where $p_i$ are
positive integers (satisfying the hypotheses of the theorem). \medskip

 Fix a valuation $v\in S$.  We shall
argue in different ways, according to the following three possibilities
for $P^v$.

{\it Case A}:  $P^v$ does not belong to the support
$|D|$ of $D$.

{\it Case B}:  $P^v$ lies in exactly one of the
irreducible components  of $|D|$, which we call ${D^v}$.

{\it Case C}:
$P^v$ lies in exactly two of the $D_i$'s, which we call ${D^v}, D^v_*$.
\medskip

Note that our assumption that no three of the $D_i$'s share a common
point implies that  no other cases may occur.\medskip

We fix an integer  $N$, sufficiently large to justify the subsequent arguments.
 We then consider
the following vector space $V=V_N$:
$$
V_N=\{\varphi\in k(X): {\rm div}(\varphi)+ND\ge 0\}.
$$
Recall that we  are assuming that each $D_i$ is defined over $k$, and in
particular we
may apply the results of the previous section.  Since $X$ is nonsingular,
whence
normal, each function in $V$ is regular on $X$ (by [H, Pro. 6.3A]).
Equivalently,
$V\subset k[X]$, i.e. every function in $V$ is a polynomial in the affine
coordinates. Let then  $\varphi_1,\ldots ,\varphi_d$ be a basis for $V$
over $k$. (For
large enough $N$, we may assume $d\ge 2$.) By the above observation,
$\varphi_j\in k[X]$, so on multiplying all  the $\varphi_j$ by a suitable
positive
integer, we may assume that all the values $\varphi_j(P_i)$ lie in $\O_S$.
\medskip

For $v\in S$, we shall construct suitable $k$-linear forms 
$L_{1v},\ldots ,L_{dv}$
in  $\varphi_1,\ldots ,\varphi_d$, linearly independent. Our aim is
to ensure that the product $\prod_{j=1}^d|L_{jv}(P_i)|_v$ is sufficiently
small with respect to the ``local height" of  the point
$(\varphi_1(P_i),\ldots ,\varphi_d(P_i))$.

More precisely,  our first aim will be to show that, for a positive number
$\mu_v$ and for all the points in a suitable infinite subsequence of 
$\{P_i\}$,  we have
$$
\prod_{j=1}^d|L_{jv}(P_i)|_v\ll
\left(\max_j(|\varphi_j(P_i)|_v)\right)^{-\mu_v},\eqno(3.1)
$$
where the implied constant does not depend on $i$.\medskip

During this construction, where $v$ is supposed to be fixed, we
shall sometimes omit the reference to it, in order to ease the
notation.\medskip

In {\bf Case A}, we simply choose $L_{jv}=\varphi_j$. Since now  all the
functions $\varphi_j$ are regular at $P^v$, they are bounded on the whole
sequence $P_i$. Therefore
$$
\prod_{j=1}^d|L_{jv}(P_i)|_v\ll
\left(\max_j(|\varphi_j(P_i)|_v)\right)^{-1},
$$
where the implied constant does not depend on $i$, and so (3.1) holds with
$\mu_v=1$.
(Note that since the constant function $1$ lies in $V$, not all the
$\varphi_j$ can vanish at $P_i$.)\medskip

We now consider {\bf Case B}, namely the sequence $\{P_i\}$ converges
$v$-adically to a point $P^v$ lying in ${D^v}$ but in no other of the
divisors $D_j$.  Since $\tilde X$ is nonsingular, we may choose,
once and for all, a local equation $t_v=0$ at $P^v$ for the divisor 
${D^v}$,  where $t_v$ is a suitable rational function on $X$. \medskip

We define a filtration of $V=V_N$ by putting
$$
W_j:=\{\varphi\in V|\ \ord_{{D^v}}(\varphi)\ge
j-1-Np^v\}, \qquad j=1,2,\ldots .\eqno(3.2)
$$
Here we put $p^v=p_i$, if ${D^v}$ is the divisor $D_i$.  Observe that in
fact we have a filtration, since $V=W_1\supset W_2\supset \ldots$, where
eventually
 $W_j=\{0\}$. Starting then from the last nonzero $W_j$, we pick a basis of it
and complete it successively to bases of the previous  spaces of the
filtration. In this way we shall eventually find a basis
$\{\psi_1,\ldots ,\psi_d\}$ of $V$ containing a basis of each given
$W_j$.

In particular, this basis contains exactly $\dim (W_j/W_{j+1})$  elements
in the set $W_j\setminus W_{j+1}$;  the order at ${D^v}$ of every such
element  is precisely $j-1-Np^v$. Hence
$$
\sum_{j=1}^d\ord_{{D^v}}(\psi_j)=\sum_{j\ge 1}(j-1-Np^v)\dim(W_j/W_{j+1}).
\eqno(3.3)
$$
Our next task is to obtain a lower bound for the right side.
To do this it will be convenient to state separately a little
combinatorial lemma.\medskip

\noindent{\bf Lemma 3.1.}\enspace {\it Let $d,U_1,\ldots ,U_h\ge 0$ and
let  $R$ be an integer $\le h$ such that
$\sum_{j=1}^RU_j\le d$. Suppose further that the real numbers
$x_1,\ldots ,x_h$ satisfy $0\le x_j\le U_j$ and $\sum_{j=1}^hx_j=d$.
Then $\sum_{j=1}^hjx_j\ge \sum_{j=1}^RjU_j$.}\medskip

\noindent{\it Proof.} We have
$$
\eqalign{\sum_{j=1}^RjU_j+\sum_{j=1}^h(R+1-j)x_j &\le
\sum_{j=1}^RjU_j+\sum_{j=1}^R(R+1-j)x_j\cr  & \le
\sum_{j=1}^RjU_j+\sum_{j=1}^R(R+1-j)U_j=(R+1)\sum_{j=1}^RU_j.}
$$
But $\sum_{j=1}^h(R+1-j)x_j=(R+1)d-\sum_{j=1}^hjx_j$, whence
$\sum_{j=1}^hjx_j\ge \sum_{j=1}^RjU_j+(R+1)(d-\sum_{j=1}^RU_j)$
and the result follows since $d-\sum_{j=1}^RU_j\ge 0$.///\medskip

We shall apply the lemma, taking $x_j:=\dim (W_j/W_{j+1})$ and
defining $h$ to be the number of nonzero $W_j$. Observe that
$\sum_{j=1}^hx_j=\dim V=d$, consistently with our previous notation.
Recall from the previous section (Lemma 2.1) that, for $D$ as in the
statement of the theorem,
 $$
d={N^2D^2\over 2} + O(N),\eqno(3.4)
$$
where the implied constant  depends only on the surface $\tilde X$ and on
the
divisor $D$.

Further, let us define $U_j=1+N(D.{D^v})-j{D^v}^2$ for $j=1,\ldots ,h$. Note
that, by Lemma 2.3, $0\le x_j\le U_j$ for $j=1,\ldots ,h$.\medskip

Let $\xi$ denote the minimal positive solution of the equation
$$
{D^v}^2\xi^2-2(D.{D^v})\xi+D^2=0,
$$
so $\xi=\xi_i$ if $D^v=D_i$. Note that by Lemma 2.4  the solutions of this equation are real, and they
cannot all be
 $\le 0$ because both $D^2$ and $D.D^v$ are positive (since $D$ is
ample). We also deduce that
$$
(D.{D^v})\ge \xi {D^v}^2.
$$
In fact, this is clear if ${D^v}^2\le 0$. Otherwise both roots must be
positive, with sum $2{(D.D^v)\over {D^v}^2}$; and the assertion again
follows since
$\xi$ is the minimal root.\medskip

We now choose $\lambda$ to be positive, $<\xi$ and such that
$$
{\lambda^2(D.{D^v})\over 2}-{\lambda^3{D^v}^2\over 3}-{D^2p^v\over
2}>0.\eqno(3.5)
$$
This will be possible by continuity, in view of the assumption (ii) of the
theorem, applied with $\xi_i=\xi$. 
In fact, by assumption we have $2D^2\xi>(D.D^v)\xi^2+3D^2p^v$.  

Now, using the equation for $\xi$ we see
that $2D^2\xi - (D.D^v)\xi^2=
3(D.D^v)\xi^2-2{D^v}^2\xi^3$. Therefore the previous inequality yields
$3(D.D^v)\xi^2-2{D^v}^2\xi^3-3D^2p^v>0$. So (3.5) will be true for all
$\lambda$ sufficiently near to $\xi$.\medskip

Also, since $\lambda <\xi$ we have, by definition of $\xi$,
$$
(D.{D^v})\lambda-{{D^v}^2\lambda^2\over 2}<{D^2\over 2}.\eqno(3.6)
$$

We shall apply Lemma 3.1, defining $R=[\lambda N]$. We first verify that 
$\sum_{j=1}^RU_j\le d$ for large enough $N$. In fact, we have
$$
\eqalign{\sum_{j=1}^RU_j & =RN(D.{D^v})-{R^2{D^v}^2\over 2}+O(R+N)\cr &
\le N^2\left((D.{D^v})\lambda-{{D^v}^2\lambda^2\over 2}\right)+O(N)}
$$
and the conclusion follows from (3.4), since by $(3.6)$ the number into
brackets is $<D^2/2$.

Observe that, since $0\le (D.{D^v})-\xi {D^v}^2\le 
(D.{D^v})-\lambda {D^v}^2$, we
have $U_j>0$ for $j\le R$, provided $N$ is large enough. Thus, if we had
$R>h$, the sum $\sum_{j=1}^RU_j$ would be strictly larger than
$\sum_{j=1}^hx_j=d$, a contradiction which proves that $R\le h$.

We may thus apply Lemma 3.1, which yields
$$
\sum_{j=1}^hjx_j\ge
\sum_{j=1}^RjU_j=\sum_{j=1}^Rj(1+N(D.{D^v})-j{D^v}^2).
$$
The right side is $N^3\left({\lambda^2(D.{D^v})\over 2}-
{\lambda^3{D^v}^2\over 3}+O(1/N)\right)$, so   we obtain from $\sum x_j=d$,
$$
\eqalign{N^{-3}\sum_{j=1}^h(j-1-Np^v)x_j &\ge
N^{-3}\left(\sum_{j=1}^hjx_j-(Np^v+1)d\right)\cr  &\ge
{\lambda^2(D.{D^v})\over
2}-{\lambda^3{D^v}^2\over 3}-{D^2p^v\over 2}+O(1/N).} 
$$
By (3.5) the right side will be positive for large $N$; together with
(3.3) this proves that, if $N$  has been chosen sufficiently large,
$$
\sum_{j=1}^d\ord_{{D^v}}(\psi_j)>0. \eqno(3.7)
$$

Now, the functions $\psi_j$ may be expressed as linear forms in the
$\varphi_\ell$. We then put $L_{jv}=\psi_j$. We have
$$
L_{jv}=t_v^{\ord_{{D^v}}(\psi_j)}\rho_{jv},
$$
where $\rho_{jv}$ are rational functions on $\tilde X$, regular at
$P^v$. In particular, the values $\rho_{jv}(P_i)$ are defined for large
$i$ and are $v$-adically bounded as $P_i$ varies. Hence
$$
\prod_{j=1}^d|L_{jv}(P_i)|_v
\ll |t_v(P_i)|_v^{\sum_{j=1}^d\ord_{{D^v}}(\psi_j)}.
$$
By a similar argument, we have
$$
\max_j|\varphi_j(P_i)|_v\ll  |t_v(P_i)|_v^{-Np_v}.
$$
Both displayed formulas make sense for all but a finite number of the
points $P_i$, which we tacitly exclude. Then, the implied constants do
not depend on $i$.

From these inequalities we finally obtain
$$
\prod_{j=1}^d|L_{jv}(P_i)|_v\ll
\left(\max_j|\varphi_j(P_i)|_v\right)^{-\mu_v},
$$
for some positive $\mu_v$ independent of $i$; therefore we have shown (3.1)
in this
case. This concludes our discussion of Case B.\medskip

        %%%%%%%%%%%%%%%%%%%%%%

We finally treat {\bf Case C}, namely the sequence $\{P_i\}$ converges
$v$-adically to a point $P^v\in {D^v}\cap D^v_*$, where ${D^v},D^v_*$ are two
distinct divisors in the set $\{D_1,\ldots ,D_r\}$. Similarly to the
above, we denote by $p^v$, $p^v_*$ the corresponding coefficients in
$D$.  

By assumption, $P^v$ cannot belong to
a third divisor in our set; let us choose two local equations
$t_v=0$ and $t_v^*=0$ for $D^v,D_*^v$ respectively. Here $t_v,t_v^*$
are regular functions, vanishing at $P^v$; also, since $D_v,D_v^*$
are distinct and irreducible, $t_v$ and $t_v^*$ are coprime in the local
ring of $\tilde{X}$ at $P^v$. 
 \medskip

We shall now consider {\it two} filtrations on the vector space $V=V_N$,
namely we put
$$ 
\eqalign{
W_j &:=\{\varphi\in V| \ord_{D^v}(\varphi)\geq j-1-Np^v\},\cr
W_j^* &:=\{\varphi\in V| \ord_{D_*^v}(\varphi)\geq j-1-Np_*^v\}.
}
$$
The following lemma from linear algebra will be used to construct 
a suitable basis for $V$. 
\medskip

\noindent {\bf Lemma 3.2.} {\it Let $V$ be  
vector space of finite dimension $d$ over a field $k$. Let 
$V=W_1\supset W_2\supset\ldots\supset W_h$, 
$V=W_1^*\supset W_2^*\supset\ldots\supset W_{h^*}$ 
be  two filtrations on $V$. 
There exists a basis $\psi_1,\ldots,\psi_d$ of $V$
which contains a basis of each $W_j$ and each $W_j^*$.
}
\medskip

\noindent {\it Proof}. We argue by induction on $d$, the case 
$d=1$ being clear. Then we can certainly suppose
(by refining the first filtration)  that $W_2$
is a hyperplane in $V$. Put $W_i^\prime:=W_i^*\cap W_2$.
By the inductive hypothesis there exists a basis $\psi_1,\ldots,\psi_{d-1}$
of $W_2$ containing basis of both $W_3,\ldots,W_h$ and 
$W_1^\prime,\ldots,W_{h^*}^\prime$. If all the $W_i^*$
for $i=2,\ldots,h^*$  are contained in $W_2$,
then $W_i^\prime=W_i^*$ for all $i>1$; in this case we just complete 
$\{\psi_1,\ldots,\psi_{d-1}\}$ to any basis of $V$ and we are done. 
Otherwise, let $l$ be the minimum index with $W_l^*\not\subset W_2$;
in this case let $\psi_d$ be any element in $W_l^*\setminus W_2$. 
We claim that the basis $\{\psi_1,\ldots,\psi_d\}$ of $V$ has the required
 property. Plainly it contains a basis of every $W_1,\ldots,W_h$. Let
$i$ be an index in $\{1,\ldots,h^*\}$; we shall prove that the set $\{\psi_1,\ldots,
\psi_d\}$ contains a basis of $W_i$. This is true by construction if $i>l$,
because in this case $W_i^*=W_i^\prime$; if $i\leq l$, then the set
$\{\psi_1,\ldots, \psi_d\}$ contains the element $\psi_d\in W_l^*\subset W_i^*$
and it contains a basis for $W_i^\prime$, which is a hyperplane in $W_i^*$; 
hence it contains a basis of $W_i^*$.

\medskip

Now, let $\psi_1,\ldots,\psi_d$ be a basis as in Lemma 3.2.
 Again, we define the linear forms $L_{jv}$ in the $\varphi_\ell$
to satisfy $L_{jv}=\psi_j$. In analogy with Case B, we may write
$$
L_{jv}=t_v^{\ord_{{D^v}}\psi_j}{t_v^*}^{\ord_{D^v_*}\psi_j}\rho_{jv}
$$
where the $\rho_{jv}\in k(X)$ are regular at $P^v$; so, as before, their
values at the $P_i$ are defined for large $i$ and $v$-adically bounded as
$i\rightarrow\infty$. Here we have used the fact that $P^v$
is a smooth point, so the corresponding local ring is a unique
factorization domain; in particular if a regular function is divisible 
both by a power of $t_v$ and  a power of $t_v^*$ (which are coprime), 
it is divisible by their product.
 
Then we have
$$
\prod_{j=1}^d|L_{jv}(P_i)|_v\ll
|t_v(P_i)|_v^{(\sum_{j=1}^d\ord_{{D^v}}\psi_j)+(\sum_{j=1}^d\ord_{{D_*^v}}\psi_j)}
$$
where the implied constant does not depend on $i$.\medskip

Again, from the assumption (ii) applied to $D^v$ and $D_*^v$, 
the same argument as in Case B gives the analogue of (3.7), both
for $\sum_{j=1}^d\ord_{{D^v}}\psi_j$ and for 
$\sum_{j=1}^d\ord_{{D_*^v}}\psi_j$. Hence, as before, we deduce (3.1).

In conclusion,   we have proved that (3.1) holds for all
$v\in S$, for suitable choices of $\mu_v >0$. 
Also, the function constantly equal to $1$ lies in $V$, so is a linear combination of
the $\varphi_j$, so $\max |\varphi_j(P_i)|_v\gg 1$. 
Thus, letting $\mu :=\min_{v\in S}\mu_v>0$, we may
write 
$$
\prod_{j=1}^d|L_{jv}(P_i)|_v\ll
\left(\max_j|\varphi_j(P_i)|_v\right)^{-\mu},\qquad v\in S. 
$$

Our theorem will now follow by a straightforward application of the Subspace
Theorem. We recall for the reader's convenience the version we are going to
apply, equivalent to the  statement  in [S, Thm. 1$D'$, p. 178].\medskip

\noindent{\bf Subspace Theorem.}\enspace {\it For an integer $d\ge 2$ and
$v\in S$, let $L_{1v},\ldots ,L_{dv}$ 
be independent linear forms in $X_1,\ldots ,X_d$ with coefficients in
$k$, and let $\epsilon >0$. Then the solutions $(x_1,\ldots ,x_d)\in {\O}_S^d$
of the inequality
$$
\prod_{v\in S}\prod_{j=1}^d|L_{jv}(x_1,\ldots ,x_d)|_v\le
H^{-\epsilon}(x_1:\ldots :x_d)
$$
lie in the union of finitely many proper linear subspaces of $k^d$.}\medskip

We apply this theorem by putting $(x_1,\ldots ,x_d)=(\varphi_1(P_i),
\ldots,\varphi_d(P_i))$. We may assume that $H(x_1:\ldots :x_d)$ tends to infinity as
$i\rightarrow\infty$, for otherwise the projective points $(x_1:\ldots :x_d)$
would all lie in a finite set, whence the nonconstant function
$\varphi_1/\varphi_2$ would be constant, equal say to $c$, 
on an infinite subsequence of the $P_i$. In this case
the theorem follows, since infinitely many points would then lie on the curve
defined on $X$ by $\varphi_1-c\varphi_2=0$.

But then for large $i$  the $(x_1,\ldots ,x_d)$ satisfy the inequality in the
statement of the Subspace Theorem, by taking for example 
$\epsilon =\mu /2$. We may then conclude that some nontrivial linear relation 
$c_1\varphi_1(P_i)+\ldots +c_d\varphi_d(P_i)=0$, with fixed coefficients
 $c_1,\ldots ,c_d$, holds on an infinite subsequence of the $P_i$. 
Again, the theorem follows since the $\varphi_j$
are linearly independent.///\bigskip

\noindent{\it Proof of Theorem 1.} We let $p_1,\ldots ,p_r$ be positive
integers as in the statement,
 namely there exists a positive integer $c$ such that
$p_ip_j(D_i.D_j)=c$ for $1\le i,j\le r$.
 We have only to check that the assumptions
(i), (ii) for the Main Theorem are verified with this choice for the $p_i$.

Assumption (i)  actually appears also in the present theorem.
To verify (ii) note that
$$
(D.D_i)={cr\over p_i},\qquad D^2=r^2c,\qquad D_i^2={c\over p_i^2},
$$
and it follows that $\xi_i=rp_i$.  Hence inequality (ii) amounts to
$2r^3cp_i>r^3cp_i+3r^2cp_i$ which is equivalent to $r\geq 4$. 
This concludes the proof.///\bigskip

\noindent{\it Proof of Corollary 1.} We start with a few reductions. First,
 by Siegel's Theorem we may assume that, given a number field $k'$, only
finitely many of the points in question are defined over $k'$.
Next, note that we may
plainly enlarge $S$  without affecting the conclusion and we now prove that
also $k$ may be enlarged; namely, it suffices to show in (ii) that {\it all
but finitely many quadratic integral points over $k$ are mapped to
$\Pr^1(k')$ by one at least of finitely many rational maps 
$\psi\in k'(\tilde C)$ of degree 2}, where $k'$ is a finite extension of $k$.

To prove this claim, assume the last  statement; we shall deduce the
Corollary from it.  Conclusion (i) of the corollary 
remains unaltered. We now show (ii); take one of
the maps $\psi$ as in the assumed conclusion. 
We may assume that it sends to $\Pr^1(k')$  the
quadratic $S$-integral points in an infinite set $\Sigma$. Note that  the
coordinate functions $X_i$ in $k[C]$, $i=1,\ldots ,m$,  satisfy by assumption
quadratic equations $X_i^2+a_iX_i+b_i=0$, where $a_i,b_i$ are rational
functions of $\psi$; by enlarging $k'$, we may then assume that 
$a_i,b_i\in k'(\psi)$. By adding new coordinates  expressed as linear 
combinations  of the original ones,  if necessary,  the  equations show that 
 $k'(C)$ has degree $\le 2$ over $k'(a_1,b_1,\ldots,a_m,b_m)$. 
This last field is contained in $k'(\psi)$, and $[k'(C):k'(\psi)]= 2$ by
assumption;  so  $k'(\psi)=k'(a_1,b_1,\ldots,a_m,b_m)$.

By the opening remark only finitely many of the points in $\Sigma$ can be
defined over $k'$;  in the sequel we tacitly disregard these points. 
By taking suitable linear combinations (over $k$) of the coordinates, 
 we may then assume that for all points
$P\in \Sigma$ and all $i=1,\ldots ,m$, $X_i(P)\not\in k'$.  
Evaluating the equations at $P\in\Sigma$ we obtain 
$X_i(P)^2+a_i(P)X_i(P)+b_i(P)=0$. 
Note that both $a_i(P), b_i(P)$ lie in $k'$, 
since we are assuming that $\psi$ sends $\Sigma$ in $k'$.
Therefore the same equations hold by replacing $X_i(P)$ 
with its conjugate over $k$:
in fact we are assuming that $X_i(P)$ are quadratic over $k$, 
but do not lie in $k'$,
and this implies that $X_i(P)$ are of exact degree $2$ over $k'$. 
But then we see that $a_i(P),b_i(P)$ actually lie in $k$. 
Consider  the field $L=k(a_1,b_1,\ldots ,a_m,b_m)$. 
Since $L\subset k'(\psi)$, we see that $L$ is the function field of a
curve over $k$, possibly reducible over $k'$.   This curve however has the
infinitely many $k$-rational points obtained by evaluating the 
$a_i,b_i$ at $P$, for $P\in\Sigma$. 
Therefore the given curve is absolutely irreducible and of genus zero
and now the existence of $k$-rational points gives $L=k(\varphi)$ for a certain
function $\varphi\in k'(\psi)$. Since $a_i(P),b_i(P)\in k$, we have
$\varphi(P)\in k$ for $P\in \Sigma$.  
 Now, $C$ is absolutely irreducible, so $k$ is algebraically
closed in $k(C)$. Therefore $[k(C):k(\varphi)]=[k'(C):k'(\varphi)]=2$, since
$k'(\varphi)=k'(\psi)$. Therefore the function $\varphi$ may be used instead of
$\psi$ to send to $\Pr^1(k)$ (rather than $\Pr^1(k')$) the points in $\Sigma$.
\medskip

We continue by observing  that the integral points on $C$
lift to integral points of a normalization, at the cost of enlarging $k$
and $S$. Therefore, in view of what has just been shown,  we may assume that
$\tilde C$ is nonsingular.

We shall then apply Theorem 1 to the surface  $\tilde X =\tilde C^{(2)}$
defined
as  the symmetric product of $\tilde C$ with itself.
(We recall from [Se2, III.14] that $\tilde X$ is in fact smooth.)
Then we have a projection map 
$\pi:\tilde C\times\tilde C\rightarrow \tilde X$ of
 degree $2$.

We let $D_i$,
$i=1,\ldots ,r$, be the image in $\tilde X$ under $\pi$ of the divisor
$A_i\times
\tilde C\subset \tilde C\times \tilde C$.

That the $D_i$ intersect transversely, and that no three of them share a
common point follows from the corresponding fact on $\tilde C\times\tilde
C$.   Also,
note that  each $D_i$ is
ample on $\tilde X$, as follows e.g. from the Nakai-Moishezon criterion.  A
fortiori,
we have that $D_1+\ldots +D_r$ is ample. Define $X:=\tilde X\setminus
(D_1\cup\ldots\cup D_r)$; then $X$ is affine and we may fix some affine
embedding.
(That the
symmetric power of an affine variety is affine follows also from a well-known
result on quotients of a variety by a finite group of automorphisms; see for
instance [Bo, Prop. 6.15].)

Note that $\pi$ restricts to a morphism from $C\times C$ in $X$.

\medskip

Let now $\{P_i\}$ be a sequence of $S$-integral points on $C$,   such that $P_i$ is
defined over a quadratic extension $k_i$ of $k$. Letting $P_i'\in C(k_i)$
be the
point conjugate to $P_i$ over $k$, we define $Q_i:=(P_i,P_i')\in C\times C$
and
$R_i:=\pi(Q_i)\in X(k_i)$.

Observe that   $R_i\in X(k)$. In fact, for any function $\varphi\in k(X)$,
we have
that $\varphi^*=\varphi\circ\pi$ is a symmetric rational function on
$C\times C$
(that is, invariant under the natural involution of $C\times C$). Therefore
$\varphi(R_i)=\varphi^*(P_i,P_i')=\varphi^*(P_i',P_i)$. This immediately
implies that
$\varphi(R_i)$ is fixed by the Galois group $Gal(\bar k/k)$, proving the claim.

Further, we note that for any $\varphi\in k[X]$, there exists a positive
integer
$m=m_\varphi$ such that all the values $m\varphi(R_i)$ are $S$-integers. In
fact, note that $\varphi^*$ is regular on $C\times C$, that is
$\varphi^*\in k[C\times
C]=k[C]\otimes_k k[C]$; this proves the contention, since for any function
$\psi\in
k[C]$, the values $\psi(P_i),\psi(P_i')$ differ from $S$-integers by a bounded denominator, as
$i$ varies.

In particular, this assertion holds taking as $\varphi$ the coordinate
functions on
$X$. So, by multiplying such coordinates by a suitable positive integer (which
amounts to apply an affine linear coordinate change on $X$) we may assume
that the
$R_i$ are integral points on $X$.\medskip

We go on by proving that the assumptions for Theorem 1 are verified in our
situation.

Note that the pull-back of $D_i$ in $\tilde C\times \tilde C$ is given by
$\pi^*(D_i)=A_i\times \tilde C +\tilde C\times A_i$. Since every two points on
a curve represent algebraically equivalent divisors, we have that the same
holds for
the $\pi^*(D_i)$. In particular, they are numerically equivalent, so the
same is true
for the $D_i$. Since we plainly have $(\pi^*(D_1).\pi^*(D_2))=2$, it
follows that
$(D_i.D_j)=1$ for all pairs $i,j$ ([B, Prop. I.8]).

In conclusion, we have verified the assumptions for Theorem 1, with $r\ge 4$
and $p_1=\ldots =p_r=c=1$.

From Theorem 1 we deduce that the $R_i$  all lie on a certain
closed curve $Y\subset X$. To prove our assertions we may now argue
separately with
each absolutely irreducible component of $Y$. Therefore we assume that the
$R_i$ are
contained in the absolutely irreducible curve $Y$, defined over a number field
containing $k$.  Since $Y$ contains the infinitely many points $R_i$, all
defined
over $k$, it follows that $Y$ is in fact defined itself over $k$. Also, $Y$
must
have genus zero and at most two points at infinity, because of Siegel's
Theorem.  In
the sequel we also suppose, as we may, that $Y$ is closed in $X$ and we let
$\tilde
Y$ be the closure of $Y$ in $\tilde X$ and  $\tilde Z=\pi^{-1}(\tilde Y)$,
$Z=\pi^{-1}(Y)=\tilde Z\setminus (\cup_{i=1}^r\pi^*(D_i))$.  \medskip

Assume first that $r\ge 5$. Then, since $\tilde Z$ is complete at least one
of the
natural projections on $\tilde C$ is surjective, whence $\# \left(\tilde Z\cap
(\cup_{i=1}^r\pi^*(D_i))\right) \ge 5$, and therefore  $\tilde Z\setminus
Z\ge 5$.
Hence $\#(\tilde Y\setminus Y)\ge 3$, since $\#\pi^{-1}(R)\le 2$ for every
$R\in
\tilde X$. But then Siegel's Theorem applies to $Y$ and contradicts the
fact that $Y$
has infinitely many integral points. This proves part (i).\medskip

From now on we suppose that $r=4$.

The case when $C$ is rational can be treated  directly, similarly to Example
1.3 above, even without appealing to the present methods.   By extending
the ground
field and $S$, $C$ may be realized as the plane quartic
$(X-\lambda)(X^2-1)Y=1$, where
$\lambda\in k$ is not $\pm 1$. Let $(x,y)$ be a quadratic $S$-integral
point on $C$.
Denoting the conjugation over $k$ with a dash, we have that
$(x-\lambda)(x'-\lambda)=:r$, $(x-1)(x'-1)=:s$, $(x+1)(x'+1)=:t$ are all
$S$-units in
$k$. Eliminating $x,x'$  gives
$2r-(\lambda+1)s+(\lambda-1)t=2(\lambda^2-1)\neq 0$. By
$S$-unit equation-theory, as in [S2, Thm. 2A] or [V, Thm. 2.3.1], this
yields some
vanishing subsum  for all but finitely many such relations. Say that e.g.
$t=2(\lambda +1)$, $2r=(\lambda+1)s$, the other cases being analogous. This
leads to
$x+x'=\lambda +1-{s\over 2}$, $xx'=\lambda+{s\over 2}$,  whence
$x^2-(\lambda +1-{s\over 2}) x+\lambda+{s\over 2}=0$,
i.e. $r ={2(x^2-(\lambda+1)x+\lambda)\over x-1}$. Then the map  given by
$x\mapsto {2(x^2-(\lambda+1)x+\lambda)\over x-1}$ 
satisfies the conclusion.\medskip

Suppose  now   that $C$ has positive genus and view $C$ as embedded in its
Jacobian $J$. For a generic point $R\in \tilde Y$, let $\{(P,Q),
(Q,P)\}=\pi^{-1}(R)\in \tilde Z$. Then $R\mapsto P+Q\in J$ is a
well-defined rational
map  from $\tilde Y$ to $J$. But $Y$ is  a rational curve, and it is
well-known that
then such a map has to be constant ([HSi, Ex. A74(b)]), say $P+Q=c$ for
$\pi(P,Q)=R\in \tilde Y$, where $c$ is independent of $R$.  We then have a
degree 2
regular map $\psi:\tilde C\rightarrow Y$ defined by $\psi(P)=\pi((P,c-P))$. It
now suffices to note that $\psi(P_i)=\pi((P_i,P'_i))=R_i$ is an
$S$-integral point in
$Y(k)$. \medskip

\noindent{\it Proof for the Addendum.} Let $\psi$ be one of the mentioned
maps, and
let $\{P_i\}_{i\in \N}$ be an infinite sequence of distinct quadratic
integral points
on $C$ such that $\psi(P_i)\in k$.  We have equations $X_i^2+a_iX_i+b_i=0$,
where
$a_i,b_i\in k(\psi)$. By changing coordinates linearly, we may assume, as
in the
argument at the beginning of the proof of the Corollary 1, that
 $k(C)$ is quadratic over $k(a_1,b_1,\ldots ,a_m,b_m)$ and that all the
values at the
$P_i$ of the affine coordinates  $X_1,\ldots ,X_m$ are of exact degree $2$
over $k$.
Then $a_j(P_i),b_j(P_i)$ are $S$-integers in $k$, for all $i,j$ in
question. The  rational map $\varphi:P\mapsto (a_1(P),b_1(P),\ldots
,a_m(P),b_m(P))$,
from $\tilde C$ to $\Pr^1$, sends $C$ to an affine curve $Y$ (over $k$) with
infinitely many $S$-integral points over $k$. This curve, whose affine ring is
$k[Y]=k[a_1,b_1,\ldots
,a_m,b_m]$,  can have at most two
points at infinity, by Siegel's Theorem. On the other hand, the above quadratic
equations for the coordinates imply that $k[C]$ is integral over $k[Y]$,
whence  all
of the (four) points at infinity of $C$ correspond to poles of some $a_i$
or $b_i$.
Therefore the $a_1,b_1,\ldots ,a_m,b_m$ have altogether at least the four poles
$A_1,\ldots ,A_4$ on $\tilde C$, and so they have at least the  poles
$\psi(A_1),\ldots ,\psi(A_4)$, viewed as rational functions of $\psi$. But
the above rational map $\varphi$ has degree $2$, whence $\psi$  factors through
it, namely $k(\psi)=k(Y)$. Therefore the curve $Y$ has at least
$\#\{\psi(A_1),\ldots ,\psi(A_4)\}$   points at infinity. By the above
conclusion
this cardinality is at most two, proving the first contention of the
addendum.

As to the second, say that $\psi(A_1)=\psi(A_2)=:\alpha$ and
$\psi(A_3)=\psi(A_4)=:\beta$. Then ${\psi-\alpha\over\psi-\beta}$ has divisor
\note{where $\psi -\infty$ is to be interpreted as $1$}
$(A_1)+(A_2)-(A_3)-(A_4)$,
yielding a relation of the mentioned type among the $(A_i)$.///\medskip

\bigskip

{\it The authors  thank Professors Enrico Bombieri, Barbara Fantechi, Angelo
Vistoli and Paul Vojta for several very helpful discussions and comments.}

\bigskip
{\bf References.}\bigskip

\item{[B]} - A. Beauville, {\it Surfaces alg\'ebriques complexes},
Ast\'erisque, n. 54, Soc. Math. de France, 1978.\smallskip

\item{[Bo]} - A. Borel, {\it Linear Algebraic Groups}, $2^{nd}$ ed.,
Springer Verlag
GTM 126, 1991.\smallskip

\item{[CZ]} - P. Corvaja, U. Zannier, A Subspace Theorem approach to
integral points
on curves, {\it C. R. Acad. Sci. Paris}, Ser. I {\bf 334} (2002),
267-271.\smallskip

\item{[EF]} - J.-H. Evertse, R. Ferretti, Diophantine inequalities on
projective
varieties, {\it International Math. Res. Notices}, to appear.\smallskip

\item{[F]} - G. Faltings, Diophantine Approximation on Abelian Varieties, {\it
Annals of Math.} {\bf 133} (1991), 549-576.\smallskip

\item{[FW]} - G. Faltings, G. W\"ustholz, Diophantine Approximations on
Projective
Varieties, {\it Inventiones Math.} {\bf 116} (1994), 109-138.\smallskip

\item{[H]} - R. Hartshorne, {\it Algebraic Geometry}, Springer-Verlag GTM 52,
1977.\smallskip

\item{[HSi]} - M. Hindry, J.H. Silverman, {\it Diophantine Geometry},
Springer-Verlag, 2000.\smallskip

\item{[L]} - S. Lang, {\it Number Theory III}, Encyclopoedia of Mathematical
Sciences, Vol. 60,  Springer-Verlag, 1991.\smallskip

\item{[S1]} - W.M. Schmidt, {\it Diophantine Approximation},
Springer-Verlag LNM
785.\smallskip

\item{[S2]} - W.M. Schmidt, {\it Diophantine Approximations and Diophantine
Equations}, Springer-Verlag LNM 1467, 1991.\smallskip

\item{[Se1]} - J-P. Serre, {\it Lectures on the Mordell-Weil Theorem}, Vieweg,
1989.\smallskip

\item{[Se2]} - J-P. Serre, {\it Algebraic groups and class fields},
Springer Verlag,
GTM 117, 1988.\smallskip

\item{[Si]} - J.H.  Silverman, {\it The Arithmetic of Elliptic Curves},
Springer-Verlag GTM 106, 1986.\smallskip

\item{[V]} - P. Vojta, {\it Diophantine Approximations and value distribution
theory}, Springer Verlag LNM  1239.\smallskip

\vfill

P. Corvaja\hfill U. Zannier

Dip. di Matematica e Informatica\hfill  Ist. Univ. Arch. - D.C.A.

Via delle Scienze\hfill S. Croce, 191

33100 - Udine (ITALY)\hfill  30135 Venezia (ITALY)

corvaja@dimi.uniud.it\hfill zannier@iuav.it

\end

\end